\documentclass[10pt,sumlimits,namelimits]{article}

\usepackage{amsthm}
\usepackage{amsmath}
\usepackage{amsfonts}
\usepackage{amssymb}
\usepackage{mathrsfs}
\usepackage{stmaryrd}
\usepackage{dsfont}
\usepackage{enumerate}
\usepackage{color}
\usepackage[colorlinks=true,urlcolor=blue]{hyperref}
\usepackage[applemac]{inputenc}

\newcommand{\N}{\mathbb{N}}

\newcommand{\R}{\mathbb{R}}
\newcommand{\C}{\mathbb{C}}

\renewcommand{\S}{\mathbb{S}}

\newcommand{\Dr}{\mathscr{D}}

\newcommand{\0}{\bs{0}}
\newcommand{\va}{\bs{a}}
\newcommand{\vb}{\bs{b}}
\newcommand{\vc}{\bs{c}}
\newcommand{\vp}{\bs{p}}
\newcommand{\f}{\bs{f}}
\newcommand{\g}{\bs{g}}
\newcommand{\h}{\bs{h}}

\newcommand{\E}{\bs{E}}
\newcommand{\F}{\bs{F}}
\newcommand{\G}{\bs{G}}
\renewcommand{\H}{\bs{H}}
\newcommand{\U}{\bs{U}}
\newcommand{\V}{\bs{V}}

\newcommand{\vu}{\bs{u}}
\renewcommand{\v}{\bs{v}}

\newcommand{\z}{\bs{z}}
\newcommand{\vphi}{\varphi}
\newcommand{\eps}{\varepsilon}
\newcommand{\la}{\lambda}
\newcommand{\dsp}{\displaystyle}
\newcommand{\ovl}{\overline}
\newcommand{\udl}{\underline}
\newcommand{\vlim}{\lim\limits}

\newcommand{\vmin}{\min\limits}

\newcommand{\vint}{\int\limits}
\newcommand{\vsum}{\sum\limits}

\newcommand{\tends}{\longrightarrow}

\newcommand{\wt}{\widetilde}

\newcommand{\li}{\llbracket}
\newcommand{\ri}{\rrbracket}

\newcommand{\loc}{\mathrm{loc}}

\renewcommand{\c}{\mathrm{c}}
\renewcommand{\d}{\mathrm{d}}

\newcommand{\dist}{\mathrm{dist}}

\newcommand{\M}{\mathrm{max}}

\renewcommand{\le}{\leqslant}
\renewcommand{\ge}{\geqslant}
\renewcommand{\O}{\Omega}
\renewcommand{\Re}{\mathrm{Re}}
\renewcommand{\Im}{\mathrm{Im}}
\newcommand{\bs}{\boldsymbol}
\newcommand{\vi}{\bs{\mathrm{i}}}
\newcommand{\p}{\prime}

\DeclareMathOperator{\supp}{supp}

\DeclareMathOperator{\divergence}{div}

\numberwithin{equation}{section}

\newtheorem{thm}{Theorem}[section]

\newtheorem{lem}[thm]{Lemma}

\theoremstyle{definition}
\newtheorem{rmk}[thm]{Remark}
\newtheorem{defi}[thm]{Definition}

\newenvironment{proof*}{\noindent{\bf Proof.}}{\qed}
\newenvironment{vproof}[1]{\noindent{\bf Proof #1}}{\qed}

\title{\huge \sc Self-similar solutions with compactly supported profile of some nonlinear Schrödinger equations}
\author{\sc Pascal Bégout$^*$ and Jes\'us Ildefonso D\'iaz$^\dagger$}
\date{}

\begin{document}

\maketitle

\begin{gather*}
\begin{array}{cc}
            ^*\mbox{Institut de Mathématiques de Toulouse \& TSE }	&	\;^\dagger\mbox{Instituto de Matem\'atica Interdisciplinar}		\\
                                          \mbox{Université Toulouse I Capitole }	&	\mbox{ Departamento de Matem\'atica Aplicada}			\\
                                                   \mbox{Manufacture des Tabacs }	&	\mbox{ Universidad Complutense de Madrid}				\\
                                                          \mbox{21, Allée de Brienne }	&	\mbox{ Plaza de las Ciencias, 3}						\\
                                  \mbox{31015 Toulouse Cedex 6, FRANCE }	&	\mbox{ 28040 Madrid, SPAIN}
\bigskip \\
\mbox{
{\footnotesize$^*$e-mail\:: }\htmladdnormallink{{\footnotesize\udl{\tt{Pascal.Begout@math.cnrs.fr}}}}{mailto:Pascal.Begout@math.cnrs.fr}}
&
\mbox{
{\footnotesize$^\dagger$e-mail\:: }\htmladdnormallink{{\footnotesize\udl{\tt{diaz.racefyn@insde.es}}}}{mailto:diaz.racefyn@insde.es}
}
\end{array}
\end{gather*}

\begin{abstract}
``\textit{Sharp localized}'' solutions (\textit{i.e.} with compact support for each given time $t)$ of a singular nonlinear type Schrödinger equation in the whole space $\R^N$ are constructed here under the assumption that they have a self-similar structure. It requires the assumption that the external forcing term satisfies that $\f(t,x)=t^{-(\vp-2)/2}\F(t^{-1/2}x)$ for some complex exponent $\vp$ and for some profile function $\F$ which is assumed to be with compact support in $\R^N.$ We show the existence of solutions of the form $\vu(t,x)=t^{\vp/2}\U(t^{-1/2}x),$ with a profile $\U,$ which also has compact support in $\R^N.$ The proof of the localization of the support of the profile $\U$ uses some suitable energy method applied to the stationary problem satisfied by $\U$ after some unknown transformation.
\end{abstract}

{\let\thefootnote\relax\footnotetext{$^\dagger$The research of J.I.~D\'iaz was partially supported by the project ref. MTM2011-26119 of the DGISPI (Spain) and the Research Group MOMAT (Ref. 910480) supported by UCM. He has received also support from the ITN \textit{FIRST} of the Seventh Framework Program of the European Community's (grant agreement number 238702)}}
{\let\thefootnote\relax}
{\let\thefootnote\relax}

\tableofcontents

\baselineskip .65cm

\section{Introduction and main result}
\label{intro}

This paper deals with the study of ``\textit{sharp localized}'' solutions of the nonlinear type Schr\"{o}dinger equation in the whole space $\R^N,$
\begin{gather}
\label{nls}
\vi\frac{\partial\vu}{\partial t}+\Delta\vu=\va|\vu|^{-(1-m)}\vu+\f(t,x),
\end{gather}
under the fundamental assumption $m\in (0,1)$ and for different choices of the complex coefficient $\va$. Here we use the notation of bold symbols for complex mathematics entities, $\vi^2=-1$ and $\Delta=\vsum_{j=1}^N\frac{\partial^2}{\partial x^2_j}$ for the Laplacian in the variables $x.$ 

By the term ``\textit{sharp localized solutions}'' we understand solutions which are more than merely the so called ``\textit{localized solutions}'' considered earlier by many authors. For instance, most of the ``\textit{localized type solutions}'' in the previous literature must vanish at infinity in an asymptotic way:
$|\vu(t,x)|\tends0$ as $|x|\tends\infty.$ They have been intensively studied mostly when some other structure property is added to the solution. It is the case of the special solutions which receive also other names such as \textit{standing waves, travelling waves, solitons,} etc.

Here we are interested on solutions which have a sharper decay when $|x|$ goes to infinity in the sense that we will require the support of the function
$\vu(t,\:.\:)$ to be a compact set of $\R^N,$ for any $t\ge0.$

We recall that equations of the type \eqref{nls} arise in many different contexts: Nonlinear Optics, Quantum Mechanics, Hydrodynamics,
etc., and that, for instance, in Quantum Mechanics the main interest concerns the case in which $\Re(\va)>0,$ $\Im(\va)=0$ (here and in which follows
 $\Re(\va)$ is the real part of the complex number $\va$ and $\Im(\va)$ is its imaginary part) and that in Nonlinear Optics the $t$ does not represent time but the main scalar variable which appears in the propagation of the wave guide direction (see Agrawal and Kivshar~\cite{ak}, p.7; Temam and Miranville \cite{MR2169020}, p.517). Sometimes equations of the type \eqref{nls} are named as Gross-Pitaevski{\u\i} type of equations in honor of two famous papers by those authors in 1961 (Gross \cite{MR0128907} and Pitaevski{\u\i} \cite{pit}). For some physical details and many references, we send the reader to the general presentations made in the books Ablowitz, Prinari and Trubatch \cite{MR2040621}, Cazenave \cite{MR2002047} and Sulem and Sulem
\cite{MR2000f:35139}.

In most of the papers on equations of the type \eqref{nls}, it is assumed that $m=3$ (the so called \textit{cubic case}). Nevertheless there are applications in which the general case $m>0$ is of interest. For instance, it is the case of the so called ``\textit{non-Kerr type equations}'' arising in the study of \textit{optical solitons} (see, e.g., Agrawal and Kivshar \cite{ak}, p.14 and following).

The case $m\in (0,1)$ has been studied before by other authors but under different points of view: some explicit self-similar solutions (the so called 
\textit{algebraic solitons}) can be found in Polyanin and Zaitsev \cite{MR2042347} (see also Agrawal and Kivshar \cite{ak}, p.33). We also mention here the series of interesting papers by Rosenau and co-authors (Kashdan and Rosenau \cite{PhysRevLett.101.264101}, Rosenau and Schuss \cite{MR2756172}) in which ``\textit{sharp localized}'' solutions are also considered with other type of statements and methods.

We also mention that the case $\Re(\va)>0$ (which corresponds to the dissipative case, also called defocusing or repulsive case, when
$\Im(\va)=0)$ must be well distinguished of the so called attractive problem (or also focusing case) in which it is assumed that $\Re(\va)<0$ (and
$\Im(\va)=0)$. See, e.g., Ablowitz, Prinari and Trubatch \cite{MR2040621}, Cazenave \cite{MR2002047}, Sulem and Sulem \cite{MR2000f:35139} and their references).

The case of complex potentials with certain types of singularities, i.e.~corresponding to the choice $\Im(\va)\ne0,$ has been previously considered by several authors, and arises in many different situations (see, for instance, Brezis and Kato \cite{MR80i:35135}, Carles and Gallo \cite{MR2765425},  LeMesurier \cite{MR2000k:35256}, Liskevich and Stollmann \cite{MR1828819} and the references therein).

Here we assume that the datum $\f$ is not zero and represents some other physical magnitude which may arise in the possible coupling with some different phenomenon: see the different chapters of Part IV of the book Sulem and Sulem \cite{MR2000f:35139}, the interaction phenomena between long waves and short waves (Benney \cite{MR0463715}, Dias and Figueira \cite{MR2339808}, Urrea \cite{urr} and their references), etc.

Obviously, the property of the compactness of the support of $\vu(t,\:.\:)$ requires the assumption that ``the support'' of the datum function $\f(t,\:.\:)$ is a compact set of $\R^N,$ for a.e. $t>0.$ Because of that, the qualitative property we consider in this paper can be understood as a ``finite speed of propagation property'' typical of linear wave equations. We point out that our treatment is very different than other ``propagation properties'' studied previously in the literature for Schrödinger equations which are formulated in terms of the spectrum of the solutions. See, e.g., the so called Anderson localization (Anderson \cite{MR0129917}), Jensen \cite{MR797050}, etc.

One of the main reasons of the study of ``\textit{sharp localized}'' solutions arises from the fact that, if we assume for the moment $\f\equiv\0,$ then 
\begin{gather*}
\frac{\partial}{\partial t}|\vu|^2+\divergence\bs J=2\Im(\va)|\vu|^{m+1},
\end{gather*}
where
\begin{gather*}
\bs J\stackrel{\mathrm{def}}{=}\left(\vu\ovl{\nabla\vu}-\ovl\vu\nabla\vu\right)=-2\Re(\vi\,\ovl\vu\,\nabla\vu),
\end{gather*}
$(\ovl\vu$ denotes the conjugate of the complex function $\vu)$ and so we get (at least formally) that 
\begin{gather*}
\frac12\frac{\d}{\d t}\int_{\R^N}|\vu(t,x)|^2\d x=\Im(\va)\int_{\R^N}|\vu(t,x)|^{m+1}\d x.
\end{gather*}
Notice that if $\Im(\va)\neq0$ then there is no mass conservation. For instance, this is the case studied by Carles and Gallo \cite{MR2765425} where they prove that actually the solution vanishes after a finite time, once that $m\in(0,1).$ More generally, it is easy to see that the two following conservation laws hold, once $a\in\R$ and $\f\equiv\0:$ if $\vu(t)\in\bs{H^1}(\R^N)\cap\bs{L^{m+1}}(\R^N)$ then we have the mass conservation
$\frac{\d}{\d t}\|\vu(t)\|_{\bs{L^2}(\R^N)}^2=0,$ moreover, if $\vu(t)\in\bs{H^2}(\R^N)\cap\bs{L^{2m}}(\R^N)$ then $\vu(t)\in\bs{L^{m+1}}(\R^N)$ and we have conservation of energy $\frac{\d}{\d t}E\big(\vu(t)\big)=0,$ where
\begin{gather*}
E\big(\vu(t)\big)=\frac12\|\nabla\vu(t)\|_{\bs{L^2}(\R^N)}^2+\frac{a}{m+1}\|\vu(t)\|_{\bs{L^{m+1}}(\R^N)}^{m+1}.
\end{gather*}
Indeed, in the first case, $\Delta\vu(t)\in\bs{H^{-1}}(\R^N)$ and $|\vu(t)|^{-(1-m)}\vu(t)\in\bs{L^\frac{m+1}{m}}(\R^N).$ It follows from the
equation~\eqref{nls} that $\frac{\partial\vu(t)}{\partial t}\in\bs{H^{-1}}(\R^N)+\bs{L^\frac{m+1}{m}}(\R^N)$ and since
$\left(\bs{H^1}(\R^N)\cap\bs{L^{m+1}}(\R^N)\right)^\star=\bs{H^{-1}}(\R^N)+\bs{L^\frac{m+1}{m}}(\R^N),$ it follows that we may take the duality product of equation~\eqref{nls} with $\vi\vu(t),$ from which the mass conservation follows. In the same way, since $\vu(t)\in\bs{L^2}(\R^N)\cap\bs{L^{2m}}(\R^N)$ and $0<m<1$, we get that $\vu(t)\in\bs{L^{m+1}}(\R^N).$ We also easily have that $\Delta\vu(t)\in\bs{L^2}(\R^N)$ and
$|\vu(t)|^{-(1-m)}\vu(t)\in\bs{L^2}(\R^N).$ It follows from the equation~\eqref{nls} that $\frac{\partial\vu(t)}{\partial t}\in\bs{L^2}(\R^N)$ and so we may take the duality product of equation~\eqref{nls} with $\frac{\partial\vu(t)}{\partial t},$ from which the conservation of energy follows.

Like in the pioneering study by Schrödinger, the condition $\Im(\va)=0$ implies that $|\vu|^2$ represents a probability density, and so the study of
``\textit{sharp localized solutions}'' becomes very relevant (recall the Heisenberg Uncertainty Principle). As we will show here (sequel of previous papers by the authors, Bégout and D\'iaz \cite{MR2214595,MR2876246}), if $m\in(0,1),$ under suitable conditions on the coefficient $\va$ (for instance for
$\Re(\va)>0$ and $\Im(\va)=0),$ it is possible to get some estimates on the support of solutions $\vu(t,x)$ showing that the probability $|\vu(t,x)|^2$ to localize a particle is zero outside of a compact set of $\R^N.$

The natural structure for searching self-similar solutions is based on the transformation $\la\longmapsto\bs{u_\la},$ where for $\la>0,$ $\vp\in\C$ and $\vu\in\bs{C}\big((0,\infty);\bs{L^1_\loc}(\R^N)\big),$ we define
\begin{gather}
\label{ula}
\bs{u_\la}(t,x)=\la^{-\vp}\vu(\la^2t,\la x), \; \forall t>0, \text{ for a.e. } x\in\R^N.
\end{gather}
Recall that since $\vp\in\C$ then
$\la^{\vp}\stackrel{\text{def}}{=}\bs{e}^{\vp\ln\la}=e^{\Re(\vp)\ln\la}\bs{e}^{\vi\Im(\vp)\ln\la}=\la^{\Re(\vp)}\bs{e}^{\vi\Im(\vp)\ln\la}$ and that
$|\la^{\vp}|=\la^{\Re(\vp)}.$ Our main assumption on the datum $\f$ is that
\begin{gather}
\label{fla}
\f(t,x)=\la^{-(\vp-2)}\f(\la^2t,\la x), \; \forall\la>0,
\end{gather}
for some $\vp\in\C,$ for any $t>0$ and almost every $x\in\R^N,$ or equivalently, that
\begin{gather}
\label{eqprof}
\f(t,x)=t^{\frac{\vp-2}{2}}\F\left(\frac{x}{\sqrt t}\right),
\end{gather}
for any $t>0$ and almost every $x\in\R^N,$ where $\F=\f(1).$ It is easy to build functions $\f$ satisfying \eqref{fla}. Indeed, for any given function $\F,$ we define $\f$ by \eqref{eqprof}. Then $\f(1)=\F$ and $\f$ satisfies \eqref{fla}. Finally, if we assume $\Re(\vp)=\frac2{1-m}$ then a direct calculation show that if $\vu$ is a solution to~\eqref{nls} then for any $\la>0,$ $\bs{u_\la}$ is also a solution to~\eqref{nls}, and conversely.

We easily check that if $\vu$ satisfies the invariance property $\vu=\bs{u_\la},$ for any $\la>0,$ then
\begin{gather}
\label{eqprou}
\vu(t,x)=t^\frac{\vp}{2}\U\left(\frac{x}{\sqrt t}\right),
\end{gather}
for any $t>0$ and almost every $x\in\R^N,$ where $\U=\vu(1).$ Thus, we arrive to the following notion:
\begin{defi}
\label{defselsim}
Let $0<m<1,$ let $\f\in\bs{C}\big((0,\infty);\bs{L^2_\loc}(\R^N)\big)$ satisfies~\eqref{fla} and let $\vp\in\C$ be such that $\Re(\vp)=\frac2{1-m}.$ A solution $\vu$ of~\eqref{nls} is said to be \textit{self-similar} if $\vu\in\bs{C}\big((0,\infty);\bs{L^2_\loc}(\R^N)\big)$ and if for any $\la>0,$ $\bs{u_\la}=\vu,$ where $\bs{u_\la}$ is defined by~\eqref{ula}. In this cases, $\vu(1)$ is called the \textit{profile} of $\vu$ and is denoted by $\U.$
\end{defi}
It follows from equation~\eqref{nls} and~\eqref{eqprou} that $\U$ satisfies
\begin{gather}
\label{U}
-\Delta\U + \va|\U|^{-(1-m)}\U - \frac{\vi\vp}{2}\U + \frac{\vi}{2}x.\nabla\U = -\F,
\end{gather}
in $\Dr^\p(\R^N),$ where $\F=\f(1).$ Conversely, if $\U\in\bs{L^2_\loc}(\R^N)$ verifies~\eqref{U}, in $\Dr^\p(\R^N),$ then the function $\vu$ defined by
\eqref{eqprou} belongs to $\bs{C}\big((0,\infty);\bs{L^2_\loc}(\R^N)\big)$ and is a self-similar solution to~\eqref{nls}, where $\f$ is defined by \eqref{eqprof} and satisfies \eqref{fla}. It is useful to introduce the unknown transformation
\begin{gather}
\label{gU}
\g(x)=\U(x)\bs{e}^{-\vi\frac{|x|^2}{8}}.
\end{gather}
Then for any $m\in\R,$ $\vp\in\C$ and $\U\in\bs{L^2_\loc}(\R^N),$ $\U$ is a solution to~\eqref{U} in $\Dr^\p(\R^N)$  if and only if
$\g\in\bs{L^2_\loc}(\R^N)$ is a solution to
\begin{gather}
\label{g}
-\Delta\g + \va|\g|^{-(1-m)}\g - \vi\frac{N+2\vp}{4}\g - \frac1{16}|x|^2\g = -\F\bs{e}^{-\vi\frac{|\:.\:|^2}{8}},
\end{gather}
in $\Dr^\p(\R^N).$ It will be convenient to study~\eqref{g} instead of~\eqref{U}. Indeed, formally, if we multiply \eqref{g} by $\pm\ovl\g$ or
$\pm\vi\ovl\g,$ integrate by parts and take the real part, one obtains some positive or negative quantities. But the same method applied to \eqref{U} gives (at least directly) nothing because of the term $\vi x.\nabla\U.$ 

Notice that if $\vp\in\C$ is such that $\Re(\vp)=\frac2{1-m}$ and if $\f\in\bs{C}\big((0,\infty);\bs{L^2}(\R^N)\big)$ and satisfies~\eqref{fla} with $\f(t_0)$ compactly supported for some $t_0>0,$ then it follows from~\eqref{fla} that for any $t>0,$ $\supp\f(t)$ is compact. Moreover, from~\eqref{eqprou}, if $\vu$ is a self-similar solution of \eqref{nls} and if $\supp\U$ is compact then for any $t>0,$ $\supp\vu(t)$ is compact. As a matter of fact, it is enough to have that $\vu(t_0)$ is compactly supported for some $t_0>0$ to have that $\vu$ satisfies \eqref{thmmain1} below and $\supp\vu(t)$ is compact, for any $t>0.$ Indeed, $\U=\vu(1)$ satisfies~\eqref{U} and by \eqref{eqprou}, $\supp\U$ and $\supp\vu(t)$ are compact  for any $t>0.$ Let $\g$ be defined by~\eqref{gU}. Then $\g$ is a solution compactly supported to~\eqref{g} and it follows the results of Section~\ref{eus} below that $\g\in\bs{H^2_\c}(\R^N).$ By~\eqref{gU}, we obtain that $\U\in\bs{H^2_\c}(\R^N)$ and we deduce easily from~\eqref{eqprou} that $\vu$ satisfies \eqref{thmmain1}.

The main result of this paper is the following.

\begin{thm}
\label{thmmain}
Let $0<m<1,$ let $\va\in\C$ be such that $\Im(\va)\le0.$ If $\Re(\va)\le0$ then assume further that $\Im(\va)<0.$ Let $\vp\in\C$ be such that
$\Re(\vp)=\frac2{1-m}$ and let $\f\in\bs{C}\big((0,\infty);\bs{L^2}(\R^N)\big)$ satisfying~\eqref{fla}. Assume also that $\supp\f(1)$ is compact.
\begin{enumerate}
	\item
	\label{thmmaina}
		If $\|\f(1)\|_{\bs{L^2}(\R^N)}$ is small enough then there exists a self-similar solution
		\begin{gather}
		\label{thmmain1}
			\vu\in\bs C\big((0,\infty);\bs{H^2}(\R^N)\big)\cap\bs{C^1}\big((0,\infty);\bs{H^1}(\R^N)\big)\cap\bs{C^2}\big((0,\infty);\bs{L^2}(\R^N)\big)
		\end{gather}
		to~\eqref{nls} such that for any $t>0,$ $\supp\vu(t)$ is compact. In particular, $\vu$ is a strong solution and verifies~\eqref{nls} for any $t>0$ in
		$\bs{L^2}(\R^N),$ and so almost everywhere in $\R^N.$
	\item
	\label{thmmainb}
		Let $R>0.$ For any $\eps>0,$ there exists $\delta_0=\delta_0(R,\eps,|\va|,|\vp|,N,m)>0$ satisfying the following property$:$ if
		$\supp\f(1)\subset\ovl B(0,R)$ and if $\|\f(1)\|_{\bs{L^2}(\R^N)}\le\delta_0$ then the profile $\U$ of the solution obtained above verifies
		$\supp\U\subset K(\eps)\subset\ovl B(0,R+\eps),$ where
		\begin{gather*}
			K(\eps)=\Big\{x\in\R^N;\; \exists y\in\supp\f(1) \text{ such that } |x-y|\le\eps \Big\},
		\end{gather*}
		which is compact.
	\item
	\label{thmmainc}
		Let $R_0>0.$ Assume now further that $\Re(\va)>0,$ $\Im(\va)=0$ and
		\begin{gather*}
			4\Im(\vp)+2\sqrt{4\Im^2(\vp)+2}\ge R_0^2.
		\end{gather*}
		Then the solution is unique in the set of functions $\bs{C}\big((0,\infty);\bs{L^2_\c}(\R^N)\big)$ whose profile $\V$ satisfies
		$\supp\V\subset\ovl B(0,R_0).$
\end{enumerate}
\end{thm}

In contrast with many other papers on self-similar solutions of equations dealing with exponents $m>1$ (see Cazenave and Weissler
\cite{MR99d:35149,MR99f:35185,MR1745480} and their references), in this paper we do not prescribe any initial data $\vu(0)$ to \eqref{nls} since we are only interested on any solution $\vu(t)$ by an external source $\f(t)$ compactly supported. Moreover, we point out that if
$\vu\in\bs{C}\big([0,\infty);\bs{L^q}(\R^N)\big)$ is a self-similar solution to~\eqref{nls}, for some $0<q\le\infty,$ then necessarily $\vu(0)=\0.$ Indeed, with help of~\eqref{eqprou}, we easily show that $\U\in\bs{L^q}(\R^N)$ and that for any $t>0,$
$\|\vu(t)\|_{\bs{L^q}(\R^N)}=t^{\frac{1}{1-m}+\frac{N}{2q}}\|\U\|_{\bs{L^q}(\R^N)},$ implying necessarily that $\vu(0)=\0.$ On the other hand, notice that if
$\vu\in\bs{C}\big([0,\infty);\Dr^\p(\R^N)\big)$ is a self-similar solution to~\eqref{nls} then one cannot expect to have $\vu(0)\in\bs{L^q}(\R^N),$ unless
$\vu(0)=\0.$ Indeed, we would have $\bs{u_\la}(0)=\vu(0)$ in $\bs{L^q}(\R^N)$ and for any $\la>0,$
$\|\vu(0)\|_{\bs{L^q}(\R^N)}=\la^{\frac{2}{1-m}+\frac{N}{q}}\|\vu(0)\|_{\bs{L^q}(\R^N)}$ and again we deduce that necessarily $\vu(0)=\0.$ More generally, the set of functions $\vu$ satisfying the invariance property,
\begin{gather*}
\forall\la>0, \text{ for a.e. } x\in\R^N, \; \bs{u_\la}(x)\stackrel{\mathrm{def}}{=}\la^{-\vp}\vu(\la x)=\vu(x),
\end{gather*}
and lying in $\bs{L^q}(\R^N)$ is reduced to $\0.$

In the special case of self-similar solution, the above arguments show that if $\f\equiv\0,$ $a\in\R$ and $\vu\in\bs{C}\big((0,\infty);\bs{L^2_\c}(\R^N)\big)$ then necessarily $\vu(t)=0,$ for any $t>0.$ Indeed, if $\vu\in\bs{C}\big((0,\infty);\bs{L^2_\c}(\R^N)\big)$ is a self-similar solution to~\eqref{nls} then its profile $\U$ belongs to $\bs{L^2}(\R^N)$ and $\vu\in\bs C^2((0;\infty)\times\R^N)$ (see Section~\ref{eus} below). So for any $t>0,$ we can multiply the above equation by $-\vi\ovl\vu(t),$ integrate by parts over $\R^N$ and take the real part. We then deduce the mass conservation,
$\frac{\d}{\d t}\|\vu(t)\|_{\bs{L^2}(\R^N)}^2=0,$ which yields with the above identity,
\begin{gather*}
\|\U\|_{\bs{L^2}(\R^N)}=\|\vu(t)\|_{\bs{L^2}(\R^N)}=t^{\frac{1}{1-m}+\frac{N}{4}}\|\U\|_{\bs{L^2}(\R^N)},
\end{gather*}
for any $t>0.$ Hence the result. As a matter of fact, if $\ell\in\{0,1,2\}$ and if $\vu\in\bs{C}\big((0,\infty);\bs{H^\ell}(\R^N)\big)$ is a self-similar solution
to~\eqref{nls} then one easily deduces from~\eqref{eqprou} that actually $\vlim_{t\searrow0}\|\vu(t)\|_{\bs{H^\ell}(\R^N)}=0.$

We also mention here that our treatment of sharp localized solutions has some indirect connections with the study of the ``unique continuation property''. Indeed, we are showing that this property does not hold when $m\in(0,1),$ in contrast to the case of linear and other type of nonlinear Schrödinger equations (see, e.g., Kenig, Ponce and Vega \cite{MR1980854}, Urrea \cite{urr}).

The paper is organized as follows. In the next section, we introduce some notations and give  general versions of the main results
(Theorems~\ref{thmsta} and \ref{thmG}). In Section~\ref{eus}, we recall some existence, uniqueness, \textit{a priori} bound and smoothness results of solutions to equation~\eqref{g} associated to the evolution equation~\eqref{nls}. Finally, Section~\ref{proof} is devoted to the proofs of the mentioned results, which we carry out by improving some energy methods presented in Antontsev, D\'iaz and Shmarev \cite{MR2002i:35001}.

\section{Notations and general versions of the main result}
\label{not}

Before stating our main results, we will indicate here some of the notations used throughout. For $1\le p\le\infty,$ $p^\p$ is the conjugate of $p$ defined by $\frac1p+\frac1{p^\p}=1.$ We denote by $\ovl\O$ the closure of a nonempty subset $\O\subseteq\R^N$ and by $\O^\c=\R^N\setminus\O$ its complement. We note $\omega\Subset\O$ to mean that $\ovl\omega\subset\O$ and that $\ovl\omega$ is a compact subset of $\R^N.$ Unless if specified, any function lying in a functional space $\big(\bs{L^p}(\O),$
$\bs{W^{m,p}}(\O),$ etc\big) is supposed to be a complex-valued function $\big(\bs{L^p}(\O;\C),$ $\bs{W^{m,p}}(\O;\C),$ etc\big).
For a functional space $\E\subset\bs{L^1_\loc}(\O;\C),$ we denote by $\bs{E_\c}=\big\{\f\in\E;\supp\f\Subset\O\big\}.$ For a Banach space $\E,$ we denote by $\E^\star$ its topological dual and by $\langle\: . \; , \: . \:\rangle_{\E^\star,\E}\in\R$ the $\E^\star-\E$ duality product. In particular, for any $\bs T\in\bs{L^{p^\p}}(\O)$ and $\bs{\vphi}\in\bs{L^p}(\O)$ with $1\le p<\infty,$
$\langle\bs T,\bs{\vphi}\rangle_{\bs{L^{p^\p}}(\O),\bs{L^p}(\O)}=\Re\vint_\O\bs T(x)\ovl{\bs{\vphi}(x)}\d x.$ For $x_0\in\R^N$ and $r>0,$ we denote by $B(x_0,r)$ the open ball of $\R^N$ of center $x_0$ and radius $r,$ by $\S(x_0,r)$ its boundary and by $\ovl B(x_0,r)$ its closure. As usual, we denote by $C$ auxiliary positive constants, and sometimes, for positive parameters $a_1,\ldots,a_n,$ write $C(a_1,\ldots,a_n)$ to indicate that the constant $C$ continuously depends only on $a_1,\ldots,a_n$ (this convention also holds for constants which are not denoted by ``$C$'').

\medskip
Now, we state the precise notion of solution.

\begin{defi}
\label{defsols}
Let $\O$ be a nonempty bounded open subset of $\R^N,$  let $(\va,\vb,\vc)\in\bs{\C^3},$ let $0<m\le1$ and let $\G\in\bs{L^1_\loc}(\O).$
\begin{enumerate}
   \item
   \label{def1}
    We say that $\g$ is a {\it local very weak solution} to
	\begin{gather}
		\label{eq1}
			-\Delta\g + \va|\g|^{-(1-m)}\g + \vb\g + \vc x.\nabla\g = \G,
	\end{gather}
    in $ \Dr^\p(\O),$ if $\g\in\bs{L^2_\loc}(\O)$ and if
	\begin{gather}
		\label{defsol0}
			\langle\g,-\Delta\bs\vphi\rangle_{\Dr^\p(\O),\Dr(\O)}+\langle\H(\g),\bs\vphi\rangle_{\Dr^\p(\O),\Dr(\O)}
			=\langle\G,\bs\vphi\rangle_{\Dr^\p(\O),\Dr(\O)},
	\end{gather}
    for any $\bs\vphi\in\Dr(\O),$ where
	\begin{gather}
	\label{H1}
		\H(\h)=\va|\h|^{-(1-m)}\h + \vb\h + \vc x.\nabla\h,
	\end{gather}
    for any $\h\in\bs{L^2_\loc}(\O).$ If, in addition, $\g\in\bs{L^2}(\O)$ then we say that $\g$ is a {\it global very weak solution} to~\eqref{eq1}.
   \item
   \label{def2}
    We say that $\g$ is a {\it local weak solution} to~\eqref{eq1} in $ \Dr^\p(\O),$ if $\g\in\bs{H^1_\loc}(\O)$ and if
	\begin{gather}
		\label{defsol1}
			\langle\nabla\g,\nabla\bs\vphi\rangle_{\Dr^\p(\O),\Dr(\O)}+\langle\H(\g),\bs\vphi\rangle_{\Dr^\p(\O),\Dr(\O)}
			=\langle\G,\bs\vphi\rangle_{\Dr^\p(\O),\Dr(\O)},
	\end{gather}
    for any $\bs\vphi\in\Dr(\O),$ where $\H\in\bs C\big(\bs{L^2_\loc}(\O);\Dr^\p(\O)\big)$ is defined by \eqref{H1}.
   \item
   \label{def3}
    We say that $\g$ is a {\it local weak solution} to
	\begin{gather}
		\label{eq2}
			-\Delta\g + \va|\g|^{-(1-m)}\g + \vb\g + \vc|x|^2\g = \G,
	\end{gather}
    in $\Dr^\p(\O),$ if $\g\in\bs{H^1_\loc}(\O)$ and if $\g$ satisfies \eqref{defsol1}, for any $\bs\vphi\in\Dr(\O),$ where
	\begin{gather}
	\label{H2}
		\H(\h)=\va|\h|^{-(1-m)}\h + \vb\h + \vc|x|^2\h,
	\end{gather}
    for any $\h\in\bs{H^1_\loc}(\O).$
   \item
   \label{def4}
    Assume further that $\G\in\bs{L^2}(\O).$ We say that $\g$ is a {\it global weak solution} to \eqref{eq1} and
	\begin{gather}
		\label{dir}
			 \g_{|\Gamma}=\0,
	\end{gather}
    in $\bs{L^2}(\O),$ if $\g\in\bs{H^1_0}(\O)$ and if
	\begin{gather}
		\label{defsol2}
			\langle\nabla\g,\nabla\v\rangle_{\bs{L^2}(\O),\bs{L^2}(\O)}+\langle\H(\g),\v\rangle_{\bs{L^2}(\O),\bs{L^2}(\O)}
			=\langle\G,\v\rangle_{\bs{L^2}(\O),\bs{L^2}(\O)},
	\end{gather}
    for any $\v\in\bs{H^1_0}(\O),$ where $\H\in\bs C\big(\bs{H^1}(\O);\bs{L^2}(\O)\big)$ is defined by \eqref{H1}. Note that $\Delta\g\in\bs{L^2}(\O),$ so
    that equation \eqref{eq1} makes sense in $\bs{L^2}(\O)$ and almost everywhere in $\O.$
   \item
   \label{def5}
    Assume further that $\G\in\bs{L^2}(\O).$ We say that $\g$ is a {\it global weak solution} to \eqref{eq2} and \eqref{dir}, in $\bs{L^2}(\O),$ if
    $\g\in\bs{H^1_0}(\O)$ and if $\g$ satisfies \eqref{defsol2}, for any $\v\in\bs{H^1_0}(\O),$ where $\H\in\bs C\big(\bs{L^2}(\O);\bs{L^2}(\O)\big)$ is
    defined by \eqref{H2}. Note that $\Delta\g\in\bs{L^2}(\O),$ so that equation \eqref{eq2} makes sense in $\bs{L^2}(\O)$ and almost everywhere in $\O.$
\end{enumerate}
\end{defi}

\noindent
In the above definition, $\Gamma$ denotes the boundary of $\O$ and $\bs C(\O)=\bs{C^0}(\O)$ is the space of complex-valued functions which are defined and continuous over $\O.$ Obviously, for $k\in\N,$ $\bs{C^k}(\O)$ denotes the space of complex-valued functions lying in $\bs C(\O)$ and having all derivatives of order lesser or equal than $k$ belonging to $\bs C(\O).$

\begin{rmk}
\label{rmkdefsols}
Here are some comments about Definition~\ref{defsols}.
\begin{enumerate}
 \item
  \label{rmkdefsols1}
	Note that in Definition~\ref{defsols}, any global weak solution is a local weak and a global very weak solution, and any local weak or global very
	weak solution is a local very weak solution.
 \item
  \label{rmkdefsols2}
	Assume that $\O$ has a $C^{0,1}$ boundary. Let $\g\in\bs{H^1}(\O).$ Then boundary condition $\g_{|\Gamma}=0$ makes sense in the sense of
	the trace $\bs\gamma(\g)=\0.$ Thus, it is well-known that $\g\in\bs{H^1_0}(\O)$ if and only if $\bs\gamma(\g)=\0.$ If furthermore $\O$ has a $C^1$
	boundary and if $\g\in\bs C(\ovl\O)\cap\bs{H^1_0}(\O)$ then for any $x\in\Gamma,$ $\g(x)=\0$ (Theorem~9.17, p.288, in Brezis \cite{MR2759829}).
	Finally, if $\g\not\in\bs C(\ovl\O)$ and $\O$ has not a $C^{0,1}$ boundary, the condition $\g_{|\Gamma}=\0$ does not
	make sense and, in this case, has to be understood as $\g\in\bs{H^1_0}(\O).$
 \item
  \label{rmkdefsols3}
	Let $0<m\le1$ and let $\z\in\C\setminus\{\0\}.$ Since $\left||\z|^{-(1-m)}\z\right|=|\z|^m,$ it is understood in Definition~\ref{defsols} that
	$\left||\z|^{-(1-m)}\z\right|=0$ when $\z=\0.$
\end{enumerate}
\end{rmk}

The main results of this section are the two following theorems implying, as a special case, the statement of Theorem~\ref{thmmain}.

\begin{thm}
\label{thmsta}
Let $\O\subset B(0,R)$ be a nonempty bounded open subset of $\R^N,$ let $0<m<1,$ let $(\va,\vb,\vc)\in\C^3$ be such that $\Im(\va)\le0,$
$\Im(\vb)<0$ and $\Im(\vc)\le0.$ If $\Re(\va)\le0$ then assume further that $\Im(\va)<0.$ Then there exist three positive constants $C=C(N,m),$
$L=L(R,|\va|,|\vp|,N,m)$ and $M=M(R,|\va|,|\vp|,N,m)$ satisfying the following property$:$ let $\G\in\bs{L^1_\loc}(\O),$ let $\g\in\bs{H^1_\loc}(\O)$ be any local weak solution to~\eqref{eq2}, let $x_0\in\O$ and let $\rho_0>0.$ If $\rho_0>\dist(x_0,\Gamma)$ then assume further that $\g\in\bs{H^1_0}(\O).$ Assume now that $\G_{|\O\cap B(x_0,\rho_0)}\equiv\0.$ Then $\g_{|\O\cap B(x_0,\rho_\M)}\equiv\0,$ where
\begin{multline}
 \label{thmsta1}
  \rho_\M^\nu=\left(\rho_0^\nu-CM^2\max\left\{1,\frac{1}{L^2}\right\}\max\left\{\rho_0^{\nu-1},1\right\}\right. \\
  \left.\times\vmin_{\tau\in\left(\frac{m+1}{2},1\right]}\left\{\frac{E(\rho_0)^{\gamma(\tau)}
  \max\{b(\rho_0)^{\mu(\tau)},b(\rho_0)^{\eta(\tau)}\}}{2\tau-(1+m)}\right\}\right)_+,
\end{multline}
where
\begin{gather*}
\begin{array}{lll}
E(\rho_0)=\|\nabla\g\|_{\bs{L^2}(\O\cap B(x_0,\rho_0))}^2,	&	&	b(\rho_0)=\|\g\|_{\bs{L^{m+1}}(\O\cap B(x_0,\rho_0))}^{m+1},	\medskip \\
k=2(1+m)+N(1-m),									&	&	\nu=\frac{k}{m+1}>2,
\end{array}
\end{gather*}
and where
\begin{gather*}
\gamma(\tau)=\frac{2\tau-(1+m)}{k}\in(0,1), \quad \mu(\tau)=\frac{2(1-\tau)}{k}, \quad \eta(\tau)=\frac{1-m}{1+m}-\gamma(\tau)>0.
\end{gather*}
for any $\tau\in\left(\frac{m+1}{2},1\right].$
\end{thm}

\noindent
Here and in what follows, $r_+=\max\{0,r\}$ denotes the positive part of the real number $r.$ 

\begin{rmk}
\label{rmkthmsta}
If the solution is too ``large'', it may happen that $\rho_\M=0$ and so the above result is not consistent. A sufficient condition to observe a localizing effect is that the solution is small enough, in a suitable sense. We give below a sufficient condition on the data $\va\in\C,$ $\vp\in\C$ and $\G$ to have
$\rho_\M>0.$
\end{rmk}

\begin{thm}
\label{thmG}
Let $\O\subset B(0,R)$ be a nonempty bounded open subset of $\R^N,$ let $0<m<1,$ let $(\va,\vb,\vc)\in\C^3$ be such that $\Im(\va)\le0,$
$\Im(\vb)<0$ and $\Im(\vc)\le0.$ If $\Re(\va)\le0$ then assume further that $\Im(\va)<0.$ Let $\G\in\bs{L^1_\loc}(\O),$ let $\g\in\bs{H^1_\loc}(\O)$ be any local weak solution to~\eqref{eq2}, let $x_0\in\O$ and let $\rho_1>0.$ If $\rho_1>\dist(x_0,\Gamma)$ then assume further that
$\g\in\bs{H^1_0}(\O).$ Then there exist two positive constants $E_\star>0$ and $\eps_\star>0$ satisfying the following property$:$ let
$\rho_0\in(0,\rho_1)$ and assume that $\|\nabla\g\|_{\bs{L^2}(\O\cap B(x_0,\rho_1))}^2<E_\star$ and
\begin{gather}
\label{thmG1}
\forall\rho\in(0,\rho_1), \;
\|\G\|_{\bs{L^2}(\O\cap B(x_0,\rho))}^2\le\eps_\star(\rho-\rho_0)_+^p,
\end{gather}
where $p=\frac{2(1+m)+N(1-m)}{1-m}.$ Then $\g_{|\O\cap B(x_0,\rho_0)}\equiv\0.$ In other words $($with the notation of Theorem~$\ref{thmsta}),$
$\rho_\M=\rho_0.$
 \end{thm}

\begin{rmk}
\label{rmkthmF}
We may estimate $E_\star$ and $\eps_\star$ as
\begin{gather*}
E_\star=E_\star\left(\|\g\|_{\bs{L^{m+1}}(B(x_0,\rho_1))}^{-1},\rho_1,\frac{\rho_0}{\rho_1},\frac{L}{M},N,m\right), \\
\eps_\star=\eps_\star\left(\|\g\|_{\bs{L^{m+1}}(B(x_0,\rho_1))}^{-1},\frac{\rho_0}{\rho_1},\frac{L}{M},N,m\right),
\end{gather*}
where $L>0$ and $M>0$ are given by Theorem~\ref{thmsta}. The dependence on $\frac{1}{\delta}$ means that if $\delta$ goes to $0$ then $E_\star$ and $\eps_\star$ may be very large. Note that $p=\frac{1}{\gamma(1)},$ where $\gamma$ is the function defined in
Theorem~\ref{thmsta}.
\end{rmk}

\section{Existence, uniqueness and smoothness}
\label{eus}

We recall the following results which are taken from other works by the authors (Bégout and D\'iaz~\cite{MR3315701}, Theorems~2.4, 2.6 and 2.12). Let $\O\subset B(0,R)$ be a nonempty bounded open subset of $\R^N,$ let $0<m<1$ and let $(\va,\vb,\vc)\in\C^3$  be such that $\Im(\va)\le0,$ $\Im(\vb)<0$ and $\Im(\vc)\le0.$ If $\Re(\va)\le0$ then assume further that $\Im(\va)<0.$ For any $\G\in\bs{L^2}(\O),$ there exists at least one global weak solution
$\g\in\bs{H^1_0}(\O)\cap\bs{H^2_\loc}(\O)$ to~\eqref{eq2} and~\eqref{dir}. Moreover, if $\O$ has a $C^{1,1}$ boundary then $\g\in\bs{H^2}(\O).$ Finally,
\begin{gather}
\label{new}
\|\g\|_{\bs{H^1}(\O)}\le M_0(R^2+1)\|\G\|_{\bs{L^2}(\O)},
\end{gather}
where $M_0=M_0(|\va|,|\vb|,|\vc|).$ Finally, if $\U$ belongs to $\bs{L^2_\loc}(\O)$ with $\U$ a local very weak solution to
\begin{gather*}
-\Delta\U + \va|\U|^{-(1-m)}\U + \vb\U + \vi cx.\nabla\U = \F, \text{ in } \Dr^\p(\O),
\end{gather*}
\big(with any $(\va,\vb,c)\in\C\times\C\times\R\big)$ then $\U\in\bs{H^2_\loc}(\O).$ Indeed, by the unknown transformation described at the beginning of Section~\ref{proof} below, we are brought back to the study of the smoothness of solutions to equation,
\begin{gather*}
-\Delta\g+\va|\g|^{-(1-m)}\g+\left(\vb-\vi\frac{cN}{2}\right)\g-\frac{c^2}{4}|x|^2\g=\F(x)\bs{e}^{-\vi c\frac{|x|^2}{4}}, \text{ in } \Dr^\p(\O),
\end{gather*}
for which the above smoothness result applies. Concerning the uniqueness of solutions, we have the following result.

\begin{thm}[\textbf{Uniqueness}]
\label{thmuni}
Let $\O\subseteq\R^N$ be a nonempty open subset let $0<m<1,$ let $(a,\vb,c)\in\R\times\C\times\R$ be such that $a>0,$ $\Re(\vb)\ge0$ and $c\ge0.$ Then for any $\F\in\bs{L^2}(\O),$ equation
\begin{gather*}
-\Delta\U - \vi a|\U|^{-(1-m)}\U - \vi\vb\U + \vi cx.\nabla\U = \F, \text{ in } \Dr^\p(\O),
\end{gather*}
admits at most one global very weak solution compact with support $\U\in\bs{L^2_\c}(\O).$
\end{thm}

\begin{proof*}
Let $\bs{U_1},\bs{U_2}\in\bs{L^2_\c}(\O)$ be two global very weak solutions both compactly supported to the above equation. By the results above, one has $\bs{U_1},\bs{U_2}\in\bs{H^2_\c}(\O).$ Setting $\bs{g_1}=\bs{U_1}\bs{e}^{-\vi c\frac{|\:.\:|^2}{4}}$ and
$\bs{g_2}=\bs{U_2}\bs{e}^{-\vi c\frac{|\:.\:|^2}{4}},$ a straightforward calculation shows that (see also the beginning of Section~\ref{proof} below)
$\bs{g_1},\bs{g_2}\in\bs{H^2_\c}(\O)$ satisfy
\begin{gather*}
-\Delta\g+\wt\va|\g|^{-(1-m)}\g+\wt\vb\g+\wt c V^2\g = \wt\F, \text{ in } \bs{L^2}(\O),
\end{gather*}
where $\wt\va=-\vi a,$ $\wt\vb=-\vi\left(\vb+\frac{cN}{2}\right),$ $\wt c=-\frac{c^2}{4},$ $V(x)=|x|$ and $\wt\F=\F\bs{e}^{-\vi c\frac{|\:.\:|^2}{4}}.$ Note that,
\begin{align*}
	&	\wt\va\neq0,	\quad	\Re(\wt\va)=0,														\\
	&	\Re\left(\wt\va\,\ovl{\wt\vb}\right)=\Re\left(a\left(\ovl{\vb+\frac{cN}{2}}\right)\right)=a\Re(\vb)+\frac12acN\ge0,	\\
	&	\Re\left(\wt\va\,\ovl{\wt c}\right)=\frac{ac^2}{4}\Re(\vi)=0.
\end{align*}
It follows from 1) of Theorem~2.10 in Bégout and D\'iaz \cite{MR3315701} that $\bs{g_1}=\bs{g_2}$ and hence, $\bs{U_1}=\bs{U_2}.$
\medskip
\end{proof*}

\begin{rmk}
\label{moreuni}
Notice that uniqueness for self-similar solution is relied to uniqueness for \eqref{g}. Using Theorem~2.10 in Bégout and D\'iaz \cite{MR3315701}, we can show that the uniqueness of self-similar solutions to equation~\eqref{nls} holds in the class of functions $\bs{C}\big((0,\infty);\bs{L^2_\c}(\R^N)\big)$ when, for instance, $\Re(\va)=0$ and $\Im(\va)<0$ (Theorem~\ref{thmuni}). These hypotheses are the same as in Carles and Gallo \cite{MR2765425}. We point out that it seems possible to adapt the uniqueness method of Theorem~2.10 in Bégout and D\'iaz \cite{MR3315701} to obtain other criteria of uniqueness.
\end{rmk}

\begin{rmk}
\label{rmkpoi}
In the proof of uniqueness of Theorem~\ref{thmmain}, we will use the Poincaré's inequality~\eqref{poibt}. This estimate can be improved in several ways. For instance, for any $x_0\in\R^N$ and any $R>0,$ we have
\begin{gather}
\label{rmkpoi1}
\|\vu\|_{\bs{L^2}(B(x_0,R))}\le\frac{2R}{\pi}\|\nabla\vu\|_{\bs{L^2}(B(x_0,R))},
\end{gather}
which is substantially better than~\eqref{poibt}, since $\frac2\pi<1<\sqrt2.$ Actually,~\eqref{rmkpoi1} holds for any $\vu\in\bs{H^1}\big(B(x_0,R)\big)$ such that
\begin{gather*}
\vint_{B(x_0,R)}\vu(x)\d x=\0,
\end{gather*}
and $\dfrac{\partial^2\vu}{\partial x_j\partial x_k}\in\bs{L^\infty}\big(B(x_0,R)\big),$ for any $(j,k)\in\li1,N\ri\times\li1,N\ri.$ See Payne and Weinberger \cite{MR0117419} for more details.
\end{rmk}

\section{Proofs of the localization properties}
\label{proof}

We start by pointing out that if $\O\subseteq\R^N$ is a nonempty open subset and if $0<m\le1,$ we have the following property: let
$\U\in\bs{H^1_\loc}(\O)$ be a local weak solution to
\begin{gather*}
-\Delta\U + \va|\U|^{-(1-m)}\U + \vb\U + \vi c x.\nabla\U = \F(x), \text{ in } \Dr^\p(\O),
\end{gather*}
for some $(\va,\vb,c)\in\C\times\C\times\R$ and $\F\in\bs{L^1_\loc}(\O).$ Setting $\g(x)=\U(x)\bs{e}^{-\vi c\frac{|x|^2}{4}},$ for almost every $x\in\O,$ it follows that $\g\in\bs{H^1_\loc}(\O)$ is a local weak solution to
\begin{gather*}
-\Delta\g+\va|\g|^{-(1-m)}\g+\left(\vb-\vi\frac{cN}{2}\right)\g-\frac{c^2}{4}|x|^2\g=\F(x)\bs{e}^{-\vi c\frac{|x|^2}{4}}, \text{ in } \Dr^\p(\O).
\end{gather*}
Conversely, if $\g\in\bs{H^1_\loc}(\O)$ is a local weak solution to
\begin{gather*}
-\Delta\g + \va|\g|^{-(1-m)}\g + \vb\g - c^2|x|^2\g = \G(x), \text{ in } \Dr^\p(\O),
\end{gather*}
for some $(\va,\vb,c)\in\C\times\C\times\R$ and $\G\in\bs{L^1_\loc}(\O),$ then setting $\U(x)=\g(x)\bs{e}^{\vi c\frac{|x|^2}{2}},$ for almost every $x\in\O,$ it follows that $\U\in\bs{H^1_\loc}(\O)$ is a local weak solution to
\begin{gather*}
-\Delta\U + \va|\U|^{-(1-m)}\U+(\vb+\vi cN)\U+2\vi c x.\nabla\U=\G(x)\bs{e}^{\vi c\frac{|x|^2}{2}}, \text{ in } \Dr^\p(\O).
\end{gather*}

The proof of Theorems~\ref{thmsta} and~\ref{thmG} follows the main structure of application of the energy methods introduced to the study of free boundary (see, e.g., the general presentation made in the monograph Antontsev, D\'iaz and Shmarev \cite{MR2002i:35001}). In both cases, the conclusions follow quite easily once it is obtained a general differential inequality for the local energy $E(\rho)$ of the type
\begin{gather}
\label{etoile}
E(\rho)^\alpha\le C\rho^{-\beta}E^\p(\rho)+K(\rho-\rho_0)_+^\omega,
\end{gather}
for some positive constants $C,$ $\beta$ and $\omega$ with $K=0,$ in case of Theorem~\ref{thmsta} and $K>0$ small enough, in case of
Theorem~\ref{thmG}. The key estimate which leads to desired local behaviour is that the exponent $\alpha$ arising in \eqref{etoile} satisfies that
$\alpha\in(0,1).$

Although the main steps to prove \eqref{etoile} follow the same steps already indicated in the monograph Antontsev, D\'iaz and Shmarev
\cite{MR2002i:35001}, it turns out that the concrete case of the systems of scalar equations generated by the Schrödinger operator does not fulfill the assumptions imposed in Antontsev, D\'iaz and Shmarev \cite{MR2002i:35001} for the case of systems of nonlinear equations. The extension of the method which applied to the system associated to the complex Schrödinger operator is far to be trivial and it was the main object of Bégout and D\'iaz \cite{MR2876246}. Unfortunately, the extension of the method presented in Bégout and D\'iaz \cite{MR2876246} is not enough to be applied to the fundamental equation of the present paper \big(\textit{i.e.} \eqref{g} or \eqref{eq2}\big) mainly due to the presence of the source term $-c^2|x|^2g.$ A sharper version of the energy method, also applicable to a different type of nonlinear complex Schrödinger type equations (for instance containing a Hartree-Fock type nonlocal term), was developed in Bégout and D\'iaz \cite{MR3190983}, where the applicability of the energy method was reduced to prove a certain local energy balance. Such a local balance will be proved here in the following lemma. Thanks to that, the proofs of Theorems~\ref{thmsta} and \ref{thmG} are then a corollary of Theorems~2.1 and 2.2 in Bégout and D\'iaz \cite{MR3190983}.

\begin{lem}
\label{lemest}
Let $\O\subset B(0,R)$ be a nonempty bounded open subset of $\R^N,$ let $0<m<1,$ let $(\va,\vb,\vc)\in\C^3$ be such that $\Im(\va)\le0,$
$\Im(\vb)<0$ and $\Im(\vc)\le0.$ If $\Re(\va)\le0$ then assume further that $\Im(\va)<0.$ Let $\G\in\bs{L^1_\loc}(\O)$ and let $\g\in\bs{H^1_\loc}(\O)$ be any local weak solution to~\eqref{eq2}. Then there exist two positive constants $L=L(R,|\va|,|\vb|,|\vc|)$ and $M=M(R,|\va|,|\vb|,|\vc|)$ such that for any $x_0\in\O$ and any $\rho_\star>0,$ if $\G_{|\O\cap B(x_0,\rho_\star)}\in\bs{L^2}\big(\O\cap B(x_0,\rho_\star)\big)$ then we have
\begin{multline}
\label{lemest1}
\|\nabla\g\|_{\bs{L^2}(\O\cap B(x_0,\rho))}^2+L\|\g\|_{\bs{L^{m+1}}(\O\cap B(x_0,\rho))}^{m+1}+L\|\g\|_{\bs{L^2}(\O\cap B(x_0,\rho))}^2		\\
\le M\left(\left|\int_{\O\cap \S(x_0,\rho)}\g\ovl{\nabla\g}.\frac{x-x_0}{|x-x_0|}\d\sigma\right|+\int_{\O\cap B(x_0,\rho)}|\G(x)\g(x)|\d x\right),
\end{multline}
for every $\rho\in[0,\rho_\star),$ where it is additionally assumed that $\g\in\bs{H^1_0}(\O)$ if $\rho_\star>\dist(x_0,\Gamma).$
\end{lem}

\begin{proof*}
Let $x_0\in\O$ and let $\rho_\star>0.$ Let $\sigma$ be the surface measure on a sphere and set for every $\rho\in[0,\rho_*),$
\begin{gather*}
I(\rho)=\left|\dsp\int_{\O\cap\S(x_0,\rho)}\g\ovl{\nabla\g}.\frac{x-x_0}{|x-x_0|}\d\sigma\right|, \quad J(\rho)=\dsp\int_{\O\cap B(x_0,\rho)}|\G(x)\g(x)|\d x, \\
w(\rho)=\int_{\O\cap\S(x_0,\rho)}\g\ovl{\nabla\g}.\frac{x-x_0}{|x-x_0|}\d\sigma, \quad
I_\Re(\rho)=\Re\big(w(\rho)\big), \quad I_\Im(\rho)=\Im\big(w(\rho)\big).
\end{gather*}
By taking as test function $\bs{\wt\vphi_n}(x)=\psi_n(|x-x_0|)\bs{\wt g}(x),$ where $\bs{\wt g}$ is the extension by $0$ of $\g$ on
$\O^\c\cap B(x_0,\rho_0)$ and $\psi_n$ is the cut-off function
\begin{gather*}
\forall t\in\R,\; \psi_n(t)=
 \begin{cases}
  \: 1,			& \mbox{ if } |t|\in\left[0,\rho-\frac{1}{n}\right], \medskip \\
  \: n(\rho-|t|),	& \mbox{ if } |t|\in\left(\rho-\frac{1}{n},\rho\right), \medskip \\
  \; 0,			& \mbox{ if } |t|\in[\rho,\infty),
 \end{cases}
\end{gather*}
it can be proved (see Theorem~3.1 in Bégout and D\'iaz \cite{MR3190983}) that $I,J,I_\Re,I_\Im\in C([0,\rho_*);\R)$ and, by passing to the limit as
$n\tends\infty,$ that
\begin{multline}
\label{prooflemest-2}
\|\nabla\g\|_{\bs{L^2}(\O\cap B(x_0,\rho))}^2+\Re(\va)\|\g\|_{\bs{L^{m+1}}(\O\cap B(x_0,\rho))}^{m+1}
										+\Re(\vb)\|\g\|_{\bs{L^2}(\O\cap B(x_0,\rho))}^2						\\
+\Re(\vc)\||x|\g\|_{\bs{L^2}(\O\cap B(x_0,\rho))}^2=I_\Re(\rho)+\Re\left(\:\vint_{\O\cap B(x_0,\rho)}\G(x)\ovl{\g(x)}\d x\right),
\end{multline}
\begin{multline}
\label{prooflemest-1}
\Im(\va)\|\g\|_{\bs{L^{m+1}}(\O\cap B(x_0,\rho))}^{m+1}+\Im(\vb)\|\g\|_{\bs{L^2}(\O\cap B(x_0,\rho))}^2
										+\Im(\vc)\||x|\g\|_{\bs{L^2}(\O\cap B(x_0,\rho))}^2					\\
=I_\Im(\rho)+\Im\left(\:\vint_{\O\cap B(x_0,\rho)}\G(x)\ovl{\g(x)}\d x\right),
\end{multline}
for any $\rho\in[0,\rho_\star).$ From these estimates, we obtain
\begin{multline}
\label{prooflemest1}
\left|\|\nabla\g\|_{\bs{L^2}(B(x_0,\rho))}^2+\Re(\va)\|\g\|_{\bs{L^{m+1}}(B(x_0,\rho))}^{m+1}+\Re(\vb)\|\g\|_{\bs{L^2}(B(x_0,\rho))}^2\right.	\\
\left.+\Re(\vc)\||x|\g\|_{\bs{L^2}(B(x_0,\rho))}^2\right|\le I(\rho)+J(\rho),
\end{multline}
\begin{gather}
\label{prooflemest2}
|\Im(\va)|\|\g\|_{\bs{L^{m+1}}(B(x_0,\rho))}^{m+1}+|\Im(\vb)|\|\g\|_{\bs{L^2}(B(x_0,\rho))}^2+|\Im(\vc)|\||x|\g\|_{\bs{L^2}(B(x_0,\rho))}^2\le I(\rho)+J(\rho),
\end{gather}
for any $\rho\in[0,\rho_\star).$ Let $A>1$ to be chosen later. We multiply~\eqref{prooflemest2} by $A$ and sum the result with~\eqref{prooflemest1}. This leads to,
\begin{multline}
\label{prooflemest3}
\|\nabla\g\|_{\bs{L^2}(B(x_0,\rho))}^2+A_1\|\g\|_{\bs{L^{m+1}}(B(x_0,\rho))}^{m+1}+A_2\|\g\|_{\bs{L^2}(B(x_0,\rho))}^2	\\
+\Re(\vc)\||x|\g\|_{\bs{L^2}(B(x_0,\rho))}^2\le2A\big(I(\rho)+J(\rho)\big),
\end{multline}
where
\begin{align*}
 A_1	&	=
		\begin{cases}
			\Re(\va),			&	\text{if } \Re(\va)>0,	\medskip \\
			A|\Im(\va)|-|\Re(\va)|,	&	\text{if } \Re(\va)\le0,
		\end{cases}
\\
 A_2	&	= A|\Im(\vb)|-|\Re(\vb)|.
\end{align*}
But~\eqref{prooflemest3} yields,
\begin{gather}
\label{prooflemest4}
\|\nabla\g\|_{\bs{L^2}(B(x_0,\rho))}^2+A_1\|\g\|_{\bs{L^{m+1}}(B(x_0,\rho))}^{m+1}+\big(A_2-R^2|\Re(\vc)|\big)\|\g\|_{\bs{L^2}(B(x_0,\rho))}^2
\le2A\big(I(\rho)+J(\rho)\big)
\end{gather}
We choose $A=A(R,|\va|,|\vb|,|\vc|)$ large enough to have $A|\Im(\va)|-|\Re(\va)|\ge1$ (when $\Re(\va)\le0)$ and $A_2-R^2|\Re(\vc)|\ge1.$
Then~\eqref{lemest1} comes from \eqref{prooflemest4} with $L=\min\big\{A_1,1\big\}$ and $M=2A.$ Note that $L=L(R,|\va|,|\vb|,|\vc|)$ and
$M=M(R,|\va|,|\vb|,|\vc|).$ This concludes the proof.
\medskip
\end{proof*}

\begin{rmk}
\label{rmklemest}
When $\rho_\star\le\dist(x_0,\Gamma)$ and $\G\in\bs{L^2_\loc}(\O),$ one may easily obtain \eqref{prooflemest-2}--\eqref{prooflemest-1} without the technical Theorem~3.1 in Bégout and D\'iaz \cite{MR3190983}. Indeed, it follows from Proposition~4.5 in Bégout and D\'iaz \cite{MR2876246} that
$\g\in\bs{H^2_\loc}(\O),$ so that equation~\eqref{eq2} makes sense in $\bs{L^2_\loc}(\O)$ and almost everywhere in $\O.$ Thus, if
$\rho_\star\le\dist(x_0,\Gamma)$ then $\g_{|B(x_0,\rho)}\in\bs{H^2}\big(B(x_0,\rho)\big)$ and \eqref{prooflemest-2} \big(respectively, \eqref{prooflemest-1}\big) is obtained by multiplying \eqref{eq2} by $\ovl\g$ (respectively, by $\ovl{\vi\g}),$ integrating by parts over $B(x_0,\rho)$ and taking the real part.
\end{rmk}

\begin{vproof}{of Theorem~\ref{thmmain}.}
Let $R>0.$ Let $\eps>0$ and let $\f\in\bs{C}\big((0,\infty);\bs{L^2}(\R^N)\big)$ satisfying~\eqref{fla} and $\supp\f(1)\subset\ovl B(0,R).$ Let $M_0$ be the constant in~\eqref{new}. Let $\vb=-\vi\frac{N+2\vp}{4},$ $\vc=-\frac1{16}$ and $\G=-\f(1)\bs{e}^{-\vi\frac{|\:.\:|^2}{8}}.$ Note that $\Im(\va)\le0,$
$\Im(\vb)=-\frac{N(1-m)+4}{4(1-m)}<0$ and $\Im(\vc)=0.$ In addition, if $\Re(\va)\le0$ then $\Im(\va)<0.$ It follows that the existence result of
Section~\ref{eus} applies to equation~\eqref{g}: let $\g\in\bs{H^1_0}(B(0,2R+2\eps))\cap\bs{H^2}(B(0,2R+2\eps))$ be such a solution to \eqref{g} and \eqref{dir}. We apply Theorem~\ref{thmsta} with $\rho_0=2\eps.$ By \eqref{new}, there exists $\delta_0=\delta_0(R,\eps,|\va|,|\vb|,|\vc|,N,m)>0$ such that if $\|\f(1)\|_{\bs{L^2}(\R^N)}\le\delta_0$ then $\rho_\M\ge\eps.$ Set $K=\supp\f(1)=\supp\G.$ Let $x_0\in\ovl{K(2\eps)^\c}\cap B(0,2R+2\eps).$ Let
$y\in B(x_0,2\eps)$ and let $z\in K.$ By definition of $K(2\eps),$ $\dist(\ovl{K(2\eps)^\c},K)=2\eps.$ We then have
\begin{gather*}
2\eps=\dist(\ovl{K(2\eps)^\c},K)\le|x_0-z|\le|x_0-y|+|y-z|<2\eps+|y-z|.
\end{gather*}
It follows that for any $z\in K,$ $|y-z|>0,$ so that $y\not\in K.$ This means that $B(x_0,2\eps)\cap K=\emptyset,$ for any
$x_0\in\ovl{K(2\eps)^\c}\cap B(0,2R+2\eps).$ By Theorem~\ref{thmsta} we deduce that for any $x_0\in\ovl{K(2\eps)^\c}\cap B(0,2R+2\eps),$
$\g_{\left|B(x_0,\eps)\right.}\equiv\0.$ By compactness, $\ovl{K(\eps)^\c}\cap B(0,2R+2\eps)$ may be covered by a finite number of sets
$B(x_0,\eps)\cap B(0,2R+2\eps)$ with $x_0\in\ovl{K(2\eps)^\c}.$ It follows that $\g_{\left| K(\eps)^\c\cap B(0,2R+2\eps)\right.}\equiv\0.$ This means that
$\supp\g\subset K(\eps)\subset B(0,2R+2\eps).$ We then extend $\g$ by $\0$ outside of $B(0,2R+2\eps).$ Thus, $\g\in\bs{H^2_\c}(\R^N)$ is a solution to \eqref{g} in $\R^N.$ Now, let $\U=\g\bs{e}^{\vi\frac{|\:.\:|^2}{8}}$ and let for any $t>0,$ $\vu(t)=t^\frac{\vp}{2}\U\left(\frac{\:\dot\:}{\sqrt t}\right).$ It follows that $\supp\U=\supp\g\subset K(\eps),$ $\U\in\bs{H^2_\c}(\R^N)$ and $\U$ is a solution to \eqref{U} in $\R^N.$ By \eqref{eqprou},
$\vu$ verifies \eqref{thmmain1} and is a solution to~\eqref{nls} in $(0,\infty)\times\R^N$ with $\vu(1)=\U$ compactly supported in
$K(\eps).$ By Definition~\ref{defselsim}, $\vu$ is self-similar and still by \eqref{eqprou}, $\supp\vu(t)$ is compact for any $t>0.$  Hence
Properties~\ref{thmmaina} and \ref{thmmainb}. It remains to show Property~\ref{thmmainc}. Let $R_0>0$ and assume further that $\Re(\va)>0,$
$\Im(\va)=0$ and $0<R_0^2\le4\Im(\vp)+2\sqrt{4\Im^2(\vp)+2}.$ Let $\bs{u_1},\bs{u_2}\in\bs{C}\big((0,\infty);\bs{L^2_\c}(\R^N)\big)$ be two solutions
to~\eqref{nls} whose profile $\bs{U_1},\bs{U_2}$ satisfy $\supp\U,\supp\V\subset\ovl B(0,R_0).$ By Section~\ref{eus},
$\bs{U_1},\bs{U_2}\in\bs{H^2_\c}(\R^N).$ For $j\in\{1,2\},$ let $\bs{g_j}=\bs{U_j}\bs{e}^{-\vi\frac{|\:.\:|^2}{8}}.$ It follows that $\bs{g_1}$ and $\bs{g_2}$ belong to $\bs{H^2_\c}(\R^N),$ are compactly supported in $\ovl B(0,R_0)$ and satisfy the same equation \eqref{g}. Let $\g=\bs{g_1}-\bs{g_2}$ and set for any $\h\in\bs{L^2_\c}(\R^N),$
$\H(\h)=|\h|^{-(1-m)}\h.$ It follows that,
\begin{gather*}
-\Delta\g + a\big(\H(\bs{g_1})-\H(\bs{g_2})\big) - \vi\frac{N+2\vp}{4}\g - \frac1{16}|x|^2\g = \0, \text{ a.e. in } \R^N.
\end{gather*}
Multiplying this equation by $\ovl\g,$ integrating by parts over $\R^N$ and taking the real part, we get
\begin{align*}
	&	\; \|\nabla\g\|_{\bs{L^2}}^2+a\langle\H(\bs{g_1})-\H(\bs{g_2}),\bs{g_1}-\bs{g_2}\rangle_{\bs{L^2},\bs{L^2}}
					 	 -\Re\left(\vi\frac{N+2\vp}{4}\right)\|\g\|_{\bs{L^2}}^2-\frac{1}{16}\||\:.\:|\g\|_{\bs{L^2}}^2		\\
   =	&	\; \|\nabla\g\|_{\bs{L^2}}^2+a\langle\H(\bs{g_1})-\H(\bs{g_2}),\bs{g_1}-\bs{g_2}\rangle_{\bs{L^2},\bs{L^2}}
									   +\frac12\Im(\vp)\|\g\|_{\bs{L^2}}^2-\frac{1}{16}\||\:.\:|\g\|_{\bs{L^2}}^2		\\
   =	&	\; 0,
\end{align*}
We recall the following refined Poincaré's inequality (Bégout and Torri \cite{BegTor}).
\begin{gather}
\label{poibt}
\forall\vu\in\bs{H^1_0}\big(B(0,R_0)\big), \; \|\vu\|_{\bs{L^2}(B(0,R_0))}^2\le2R_0^2\|\nabla\vu\|_{\bs{L^2}(B(0,R_0))}^2,
\end{gather}
If follows from~\eqref{poibt} and Lemma~9.1 in Bégout and D\'iaz \cite{MR2876246}, that there exists a positive constant $C$ such that,
\begin{gather*}
\left(\frac{1}{2R_0^2}+\frac12\Im(\vp)-\frac{R_0^2}{16}\right)\|\g\|_{\bs{L^2}}^2
+Ca\vint_\omega\frac{|\bs{g_1}(x)-\bs{g_2}(x)|^2}{(|\bs{g_1}(x)|+|\bs{g_2}(x)|)^{1-m}}\d x
\le0,
\end{gather*}
where $\omega=\Big\{x\in\O;|\bs{g_1}(x)|+|\bs{g_2}(x)|>0\Big\}.$ But,
\begin{gather*}
\frac{1}{2R_0^2}+\frac12\Im(\vp)-\frac{R_0^2}{16}=\frac{1}{16R_0^2}\left(-R_0^4+8\Im(\vp)R_0^2+8\right)\ge0,
\end{gather*}
when
\begin{gather*}
0\le R_0^2\le4\Im(\vp)+2\sqrt{4\Im^2(\vp)+2}.
\end{gather*}
It follows that $\bs{g_1}=\bs{g_2}$ which implies that $\bs{U_1}=\bs{U_2}$ and for any $t>0,$ $\bs{u_1}(t)=\bs{u_2}(t).$ This ends the proof.
\medskip
\end{vproof}

\baselineskip .4cm

\bibliographystyle{abbrv}
\bibliography{Paper10}

\def\cprime{$'$}
\begin{thebibliography}{10}

\bibitem{MR2040621}
M.~J. Ablowitz, B.~Prinari, and A.~D. Trubatch.
\newblock {\em Discrete and continuous nonlinear {S}chr\"odinger systems},
  volume 302 of {\em London Mathematical Society Lecture Note Series}.
\newblock Cambridge University Press, Cambridge, 2004.

\bibitem{ak}
G.~P. Agrawal and Y.~S. Kivshar.
\newblock {\em Optical {S}olitons: {F}rom {F}ibers to {P}hotonic {C}rystals}.
\newblock Academic Press, California, San Diego, 2003.

\bibitem{MR0129917}
P.~W. Anderson.
\newblock Localized magnetic states in metals.
\newblock {\em Phys. Rev. (2)}, 124:41--53, 1961.

\bibitem{MR2002i:35001}
S.~N. Antontsev, J.~I. D{\'{\i}}az, and S.~Shmarev.
\newblock {\em Energy methods for free boundary problems: Applications to
  nonlinear PDEs and fluid mechanics}.
\newblock Progress in Nonlinear Differential Equations and their Applications,
  48. Birkh\"auser Boston Inc., Boston, MA, 2002.

\bibitem{MR2214595}
P.~B{\'e}gout and J.~I. D{\'{\i}}az.
\newblock On a nonlinear {S}chr\"odinger equation with a localizing effect.
\newblock {\em C. R. Math. Acad. Sci. Paris}, 342(7):459--463, 2006.

\bibitem{MR2876246}
P.~B{\'e}gout and J.~I. D{\'{\i}}az.
\newblock Localizing estimates of the support of solutions of some nonlinear
  {S}chr\"odinger equations --- {T}he stationary case.
\newblock {\em Ann. Inst. H. Poincar\'e Anal. Non Lin\'eaire}, 29(1):35--58,
  2012.

\bibitem{MR3190983}
P.~B{\'e}gout and J.~I. D{\'{\i}}az.
\newblock A sharper energy method for the localization of the support to some
  stationary {S}chr\"odinger equations with a singular nonlinearity.
\newblock {\em Discrete Contin. Dyn. Syst.}, 34(9):3371--3382, 2014.

\bibitem{MR3315701}
P.~B{\'e}gout and J.~I. D{\'{\i}}az.
\newblock Existence of weak solutions to some stationary {S}chr\"odinger
  equations with singular nonlinearity.
\newblock {\em Rev. R. Acad. Cienc. Exactas F\'\i s. Nat. Ser. A Math. RACSAM},
  109(1):43--63, 2015.

\bibitem{BegTor}
P.~B{\'e}gout and V.~Torri.
\newblock Numerical computations of the support of solutions of some localizing
  stationary nonlinear {S}chr\"odinger equations.
\newblock In preparation.

\bibitem{MR0463715}
D.~J. Benney.
\newblock A general theory for interactions between short and long waves.
\newblock {\em Studies in Appl. Math.}, 56(1):81--94, 1976/77.

\bibitem{MR2759829}
H.~Brezis.
\newblock {\em Functional analysis, {S}obolev spaces and partial differential
  equations}.
\newblock Universitext. Springer, New York, 2011.

\bibitem{MR80i:35135}
H.~Brezis and T.~Kato.
\newblock Remarks on the {S}chr\"odinger operator with singular complex
  potentials.
\newblock {\em J. Math. Pures Appl. (9)}, 58(2):137--151, 1979.

\bibitem{MR2765425}
R.~Carles and C.~Gallo.
\newblock Finite time extinction by nonlinear damping for the {S}chr\"odinger
  equation.
\newblock {\em Comm. Partial Differential Equations}, 36(6):961--975, 2011.

\bibitem{MR2002047}
T.~Cazenave.
\newblock {\em Semilinear {S}chr\"odinger equations}, volume~10 of {\em Courant
  Lecture Notes in Mathematics}.
\newblock New York University Courant Institute of Mathematical Sciences, New
  York, 2003.

\bibitem{MR99d:35149}
T.~Cazenave and F.~B. Weissler.
\newblock Asymptotically self-similar global solutions of the nonlinear
  {S}chr\"odinger and heat equations.
\newblock {\em Math. Z.}, 228(1):83--120, 1998.

\bibitem{MR99f:35185}
T.~Cazenave and F.~B. Weissler.
\newblock More self-similar solutions of the nonlinear {S}chr\"odinger
  equation.
\newblock {\em NoDEA Nonlinear Differential Equations Appl.}, 5(3):355--365,
  1998.

\bibitem{MR1745480}
T.~Cazenave and F.~B. Weissler.
\newblock Scattering theory and self-similar solutions for the nonlinear
  {S}chr\"odinger equation.
\newblock {\em SIAM J. Math. Anal.}, 31(3):625--650 (electronic), 2000.

\bibitem{MR2339808}
J.-P. Dias and M.~Figueira.
\newblock Existence of weak solutions for a quasilinear version of {B}enney
  equations.
\newblock {\em J. Hyperbolic Differ. Equ.}, 4(3):555--563, 2007.

\bibitem{MR0128907}
E.~P. Gross.
\newblock Structure of a quantized vortex in boson systems.
\newblock {\em Nuovo Cimento (10)}, 20:454--477, 1961.

\bibitem{MR797050}
A.~Jensen.
\newblock Propagation estimates for {S}chr\"odinger-type operators.
\newblock {\em Trans. Amer. Math. Soc.}, 291(1):129--144, 1985.

\bibitem{PhysRevLett.101.264101}
E.~Kashdan and P.~Rosenau.
\newblock Compactification of nonlinear patterns and waves.
\newblock {\em Phys. Rev. Lett.}, 101(26):261602, 4, 2008.

\bibitem{MR1980854}
C.~E. Kenig, G.~Ponce, and L.~Vega.
\newblock On unique continuation for nonlinear {S}chr\"odinger equations.
\newblock {\em Comm. Pure Appl. Math.}, 56(9):1247--1262, 2003.

\bibitem{MR2000k:35256}
B.~J. LeMesurier.
\newblock Dissipation at singularities of the nonlinear {S}chr\"odinger
  equation through limits of regularisations.
\newblock {\em Phys. D}, 138(3-4):334--343, 2000.

\bibitem{MR1828819}
V.~Liskevich and P.~Stollmann.
\newblock Schr\"odinger operators with singular complex potentials as
  generators: existence and stability.
\newblock {\em Semigroup Forum}, 60(3):337--343, 2000.

\bibitem{MR0117419}
L.~E. Payne and H.~F. Weinberger.
\newblock An optimal {P}oincar\'e inequality for convex domains.
\newblock {\em Arch. Rational Mech. Anal.}, 5:286--292 (1960), 1960.

\bibitem{pit}
L.~P. Pitaevski{\u\i}.
\newblock Vortex lines in an imperfect bose gas.
\newblock {\em Soviet Physics. JETP}, 13(2):489--506, 1961.

\bibitem{MR2042347}
A.~D. Polyanin and V.~F. Zaitsev.
\newblock {\em Handbook of nonlinear partial differential equations}.
\newblock Chapman \& Hall/CRC, Boca Raton, FL, 2004.

\bibitem{MR2756172}
P.~Rosenau and Z.~Schuss.
\newblock Tempered wave functions: {S}chr\"odinger's equation within the light
  cone.
\newblock {\em Phys. Lett. A}, 375(5):891--897, 2011.

\bibitem{MR2000f:35139}
C.~Sulem and P.-L. Sulem.
\newblock {\em The nonlinear {S}chr\"odinger equation}, volume 139 of {\em
  Applied Mathematical Sciences}.
\newblock Springer-Verlag, New York, 1999.
\newblock Self-focusing and wave collapse.

\bibitem{MR2169020}
R.~Temam and A.~Miranville.
\newblock {\em Mathematical modeling in continuum mechanics}.
\newblock Cambridge University Press, Cambridge, second edition, 2005.

\bibitem{urr}
J.~J. Urrea.
\newblock On the support of solutions to the {N}{L}{S}-{K}{D}{V} system.
\newblock {\em Differential Integral Equations}, 25(7-8):611--6186, 2012.

\end{thebibliography}
\addcontentsline{toc}{section}{References}

\end{document}